\documentclass{ifacconf}

\usepackage{graphicx}      
\usepackage{natbib}        

\usepackage[latin1]{inputenc}
\usepackage{amssymb}

\usepackage{enumerate}

\usepackage{mathtools}

\def\begquo{\begin{quote}}
	\def\endquo{\end{quote}}
\def\begequarr{\begin{eqnarray}}
	\def\endequarr{\end{eqnarray}}
\def\begequarrs{\begin{eqnarray*}}
	\def\endequarrs{\end{eqnarray*}}
\def\begarr{\begin{array}}
	\def\endarr{\end{array}}
\def\begequ{\begin{equation}}
	\def\endequ{\end{equation}}
\def\lab{\label}
\def\begdes{\begin{description}}
	\def\enddes{\end{description}}
\def\begenu{\begin{enumerate}}
	\def\begite{\begin{itemize}}
		\def\endite{\end{itemize}}
	\def\endenu{\end{enumerate}}

\def\lef[{\left[\begin{array}}
	\def\rig]{\end{array}\right]}

\def\begcen{\begin{center}}
	\def\endcen{\end{center}}
\def\begdef{\begin{definition}}
	\def\enddef{\end{definition}}
\def\begsubequ{\begin{subequations}}
	\def\endsubequ{\end{subequations}}

\def\begmat#1{\begin{bmatrix}#1\end{bmatrix}}
\def\begali#1{\begin{align}{#1}\end{align}}
\def\begalis#1{\begin{align*}{#1}\end{align*}}

\def\caly{{\mathcal Y}}

\def\calb{{\mathcal B}}

\def\cala{{\mathcal A}}

\def\cald{{\mathcal D}}

\def\bfp{{\bf p}}

\def\L2e{{\cal L}_{2e}}

\def\rea{\mathbb{R}}

\def\adj{\mbox{adj}}
\def\col{\mbox{col}}

\usepackage{color}

\usepackage[prependcaption,colorinlistoftodos]{todonotes}

\begin{document}
\begin{frontmatter}

\title{State Observation of Affine-in-the-States  Systems with Unknown Time-Varying Parameters and Output Delay} 

\thanks[footnoteinfo]{This work is supported by the Russian Science Foundation under grant 22-21-00499 at ITMO University}

\author[First,fifth]{Alexey Bobtsov} 
\author[First]{Nikolay Nikolaev} 
\author[Third]{Romeo Ortega}
\author[Fourth]{Denis Efimov}
\author[First]{Olga Kozachek}

\address[First]{Department of Control Systems and Robotics, ITMO University, Kronverkskiy av. 49, Saint-Petersburg, 197101, Russia, (e-mail: bobtsov@mail.ru, nanikolaev@itmo.ru, oakozachek@mail.ru).}
\address[Third]{Departamento Acad\'{e}mico de Sistemas Digitales, ITAM, Ciudad de M\'exico, M\'{e}xico, (e-mail: romeo.ortega@itam.mx)}
\address[Fourth]{INRIA, Univ. Lille, CNRS, UMR 9189 - CRIStAL, F-59000 Lille, France, (e-mail: denis.efimov@inria.fr)}
\address[fifth]{Laboratory ``Control of Complex Systems", Institute of Problems of Mechanical Engineering, V.O., Bolshoj pr., 61, St. Petersburg, 199178, Russia}

\begin{abstract}                
In this paper we address the problem of adaptive state observation of affine-in-the-states time-varying systems with delayed measurements and {\em unknown} parameters. The development of the results proposed in the \citep{BOB_micnon_21} and in the \citep{Bob_et_all_IJC_21} is considered.  The case with known parameters has been studied by many researchers---see \citep{SANGARKRS,BOBetalaut21} and references therein---where, similarly to the approach adopted here, the system is treated as a {\em linear time-varying} system.  We show that the parameter estimation-based observer (PEBO) design proposed in \citep{ORTetalscl15,ORTetalaut20}  provides a very simple solution for the unknown parameter case. Moreover, when PEBO is combined with the dynamic regressor extension and mixing (DREM) estimation technique \citep{ARAetaltac17,ORTetalaut19}, the estimated state converges in {\em fixed-time} with extremely weak excitation assumptions.
\end{abstract}


\end{frontmatter}
%
\section{Introduction and Problem Formulation}
\lab{sec1}
%
It is common in control applications that real devices provide measurements with time-varying delays. This fact makes the problem of the state estimation for a dynamical system more complicated. The reason is the lags in output injection signal. This problem has been explored by many authors. In case of linear time invariant (LTI) systems, this issue is well understood and the observer convergence can be verified by checking the feasibility of a linear matrix inequality {\citep{FRI}}. On the other hand, for linear time-varying (LTV) systems this problem is widely open---see the literature review and references in the recent papers \citep{BOBetalaut21,RUEetalijc19,SANGARKRS}.
In \citep{BOB_micnon_21} and \citep{Bob_et_all_IJC_21} time-varying systems with unknown parameters were considered. In \citep{Bob_et_all_IJC_21} unknown parameters were assumed to be time-varying but the output signal was not delayed. In \citep{BOB_micnon_21} the case with delay in the output signal was considered. However unknown parameters were assumed to be constant. In this paper we develop both results from \citep{BOB_micnon_21} and \citep{Bob_et_all_IJC_21} under condition that the unknown parameters are time-varying and the output variable contains delay function.

It is well-known \citep{BERbook},  and briefly explained below, that affine-in-the-states time-varying systems are possible to be treated as LTV systems under some reasonable assumptions. The approach adopted in this paper is based on this fact. To be precise, considered affine-in-the-states time-varying systems are described by
\begali{
	\nonumber
	\dot{x}(t)&=\cala(u,y,t)x(t)+\cald(u,y,t) \eta(t) +\calb(u,y,t)\\
	\lab{sys}
	y(t)&=C(\varphi(t))x(\varphi(t)),
}
where $x(t) \in \rea^n$, $u(t) \in \rea^m$, $y(t) \in \rea^p$, $\eta(t) \in \rea^q$ is vector of {\em unknown} time-varying parameters, which can be represented in the following form
\begin{align}
	\label{eta}
	\eta(t)&=H(t) \xi(t)\\
		\dot{\xi}(t)&=\Gamma(t)\xi(t),
	\label{xi}
\end{align}
where $\Gamma(t)$ and $H(t)$ are known matrices ($H(t)$ is bounded), $\xi(t) \in \rea^k$, $\varphi(t)$ is a continuous {\em known} nonnegative function that defines the measurement delay, which is
	\begin{equation*}
		\lab{varphi}
		\varphi(t):=t-d(t),\;\varphi(t)\geq0,
	\end{equation*}
	where $d(t)$ is the time-varying delay that verifies the following inequality:
	\begequ
	\lab{boudt}
	0 \leq d(t) \leq d_M,
	\endequ
for unknown upper bound $d_M>0$, $C(\varphi(t))$ is bounded	and the state initial condition is $x(0)=x_0 \in \rea^n$. 
	
Following the standard procedure, we rewrite the state dynamics of \eqref{sys}	 as an LTV system
\begali{
	\dot{x}(t)&=A(t)x(t)+D(t) \eta(t) +B(t),
	\lab{ltvsys}
}
where we define the matrices
\begalis{
	A(t)&:=\cala(u(t),y(t),t),\\
	D(t)&:=\cald(u(t),y(t),t),\\
	B(t)&:=\calb(u(t),y(t),t).
}
We assume that these matrix functions are  {\em known}, continuous, bounded, and  verify the following assumptions.

\begin{assum} 
	\lab{ass1}
	The state transition matrices of the homogeneous part of the system \eqref{ltvsys} and the system \eqref{xi}, denoted $\Phi(t,t_0)$ and $\Phi_\Gamma(t,t_0)$, respectively, verify
	$$
	\| \Phi(t,t_0)\| +\|\Phi_\Gamma(t,t_0)\| \leq c_1,\;\forall t \geq t_0, \; c_1\in \rea.
	$$ 
\end{assum}

\begin{assum} 
	\lab{ass2}\em
	\begin{align*}
		\int_{t_0}^t \|\Phi(t,\tau)[D(\tau) H(\tau)\Phi_\Gamma(\tau,t_0)\xi(0)+B(\tau)]\| d\tau \leq c_2,\\ 
		\;\forall t \geq t_0, \; c_2\in \rea.
	\end{align*}
	\end{assum}

	Assumption \ref{ass1} is equivalent to {\em uniform stability} of the homogeneous part of the system \eqref{ltvsys} \citep[Theorem 6.4]{RUGbook} and the system \eqref{xi}, while Assumption \ref{ass2} is a necessary and sufficient condition for {\em bounded-input-bounded-state} stability of  \eqref{ltvsys}  \citep[Theorem 12.2]{RUGbook}.

	Solution of the equation \eqref{xi} we can find in the next form
	\begin{align}
		\xi(t)&=\Phi_\Gamma(t) \xi(0)=\Phi_\Gamma(t) \theta_\Gamma,\\
		\dot{\Phi}_\Gamma(t)&=\Gamma(t)\Phi_\Gamma(t), \; \Phi_\Gamma(0)=I_n,
	\end{align}
	where $\Phi_\Gamma(t)$ is the fundamental matrix and $I_n$ is the identity $n\times n$ matrix.
	
	Then equation \eqref{ltvsys} we can rewrite as follows
	\begin{align}
		\dot{x}(t)&=A(t)x(t)+ M(t) \theta_\Gamma  +B(t)
	\end{align} 
	where $M(t)=D(t)H(t)\Phi_\Gamma(t)$.

In this paper we design an adaptive observer

\begin{gather*}
\dot \chi(t) = F(\chi(t),u(t),y(t)),\\
(\hat x(t),\hat\eta(t))=K(\chi(t),u(t),y(t))
\end{gather*}

with $\chi(t) \in \rea^{n_\chi}$ such that all signals are bounded and {\em fixed-time convergence} (FCT) of the estimated state and parameters is ensured to their ideal values, that is,
	\begequ
	\lab{concon}
	\hat x(t) = x(t), \hat\eta(t)=\eta(t), \; \forall t \geq t_c, 
	\endequ
	for some $t_c \in (0,\infty)$ and for all $x_0 \in \rea^n,\xi(0)\in\rea^k,\chi(0) \in \rea^{n_\chi}$.

To design an adaptive observer we use two techniques - PEBO and DREM. PEBO is a new family of observers called parameter estimation based observers. PEBO is applicable for the class of systems that can be identified via two assumptions. 
The first one characterizes, via the solvability of a partial differential equation (PDE), systems for which there exists a partial change of coordinates that assigns a particular cascaded structure to the model, which allows to obtain a classical regression form containing the unknown parameter and measurable quantities only (it can be realized by using the flatness approach \citep{Fliess}). 
The second assumption is related with the ability to  estimate this unknown parameter consistently. This unknown parameter, in general, may enter in the regression form nonlinearly. Linear regression forms can be obtained via over parametrisation of the nonlinear  regression. There are many established parameter estimation algorithms available for linear regression forms and the second assumption can be replaced by the well known persistency of excitation (PE) condition \citep{ORTetalscl15}. DREM procedure consists of two stages. At the first one new regression forms are generated via the application of a dynamic operator to the data of the original regression. At the second one these new data are mixed suitably which allows to obtain the final desired regression form. After that standard parameter estimation techniques can be applied to this form. \citep{ARAetaltac17}.
%
\section{Main Result}
\lab{sec2}
%
To streamline the presentation of our observer we  recall from \citep{KRERIE} the following.

\begin{defn}
	\lab{de1}
	A bounded signal $\Delta :\rea_+\to \rea \in \rea$ is interval exciting (IE) if there exists a time $t_{\tt IE} \in (0,\infty)$ such that
	\begequ
	\lab{intexc}
	\int_0^{t_{\tt IE}} \Delta^2(s)ds \geq \rho,
	\endequ
	for some $\rho >0$.
\end{defn}

\begin{prop}
	{Let us consider the system \eqref{eta}, \eqref{xi} and \eqref{ltvsys} {satisfying Assumptions \ref{ass1} and \ref{ass2}} and $d(t)$ verifying \eqref{boudt}} with the generalized PEBO+DREM
	\begsubequ
	\lab{gpebodyn}
	\begali{
		\lab{dotxi}
		\dot{\zeta}(t)&=A(t)\zeta(t)+B(t)\\
			\lab{dotG}
		\dot{G}(t)&=A(t)G(t)+M(t)\\
		\lab{dotphi}
		\dot{\Phi}_A(t)&=A(t)\Phi_A(t),\;\Phi_A(0)=I_n\\
		\lab{doty}
		\dot Y(t) &= - \lambda Y(t) +  \lambda\Psi^\top (t)z(t)\\
		\lab{dotome}
		\dot \Omega(t) &=- \lambda \Omega(t) +  \lambda \Psi^\top(t)\Psi(t)
	}
	\endsubequ
	with  $\lambda >0$ and the gradient parameter estimator
	\begali{
		\lab{dothatthe}
		\dot {\hat \theta}(t) &=-\gamma \Delta(t) [\Delta(t) \hat{\theta}(t)-\caly(t)]
	}
	with $\gamma>0$, and the definitions
	\begsubequ
	\lab{gpebodyn1}
	\begali{
		\lab{psi}
		\Psi(t)&:=C(\varphi(t)) \left[\Phi_A(\varphi(t))-G(\varphi(t))\right]\\
		\lab{z}
		z(t)&:=C(\varphi(t))\zeta(\varphi(t))-y(t)\\
		\lab{caly}
		\caly(t) &:= \adj\{\Omega(t)\}Y(t)\\
		\Delta(t)&:=\det\{\Omega(t)\},
		\lab{del}
	}
	\endsubequ
where $\adj\{\cdot\}$ is the adjugate matrix. Define the FCT parameters estimate
$$
\theta^{FCT}(t)={1 \over 1 - w_c(t)}[\hat \theta(t) - w_c(t) \hat \theta(0)],
$$	
with
	\begali{
		\dot w(t)  &= -\gamma \Delta^2(t) w(t), \; w(0)=1,
		\lab{dotw}
	}
and $w_c(t)$ defined via the clipping function
	\begequ
	\lab{wc}
	w_c(t) = \left\{ \begin{array}{lcl} w(t) & \;\mbox{if}\; & w(t) \leq 1-\mu \\ 1-\mu & \;\mbox{if}\; & w(t) > 1-\mu, \end{array} \right.,
	\endequ
where $\mu \in (0,1)$ is a designer chosen parameter. Assume $\Delta(t)$ verifies \eqref{intexc} with
	\begin{align}
		\label{rho}
		\rho=-\frac{1}{\gamma}\ln(1-\mu).
	\end{align}
	Then,  the estimates
	\begali{
		\hat x(t) &= \zeta(t) {-} {\left[\Phi_{A}(t) \;|\; -G(t) \right] \theta^{FCT}(t)},
	\lab{ftcgpebo}	\\
				\hat \eta(t) &= H(t)\Phi_\Gamma(t){\left[0 \;|\; I_k \right] \theta^{FCT}(t)},
				\lab{ftcgpebo_eta}
	}
ensure \eqref{concon} for some $t_c \geq t_{\tt IE}$, and all signals remain bounded.
\end{prop}

\begin{pf}
{Firstly, we prove the claim of the signal boundedness. Due to Assumptions \ref{ass1} and \ref{ass2} 
we have that $x(t)$, $\zeta(t)$, $G(t)$ and the fundamental matrices $\Phi_A(t), \Phi_\Gamma(t)$ are bounded. Considering that $C(\varphi(t))$ is also bounded, this guarantees that $y(t)$, $z(t)$ and $\Psi(t)$ are bounded. Finally, positivity of $\lambda$ ensures that $Y(t)$ and $\Omega(t)$ are bounded as well. } 
	
	Secondly, we apply the PEBO technique to derive a vector linear regression equation (LRE). For this purpose, let us define the error signal 
	$$
	e(t):=\zeta(t)-x(t)+G(t)\theta_\Gamma,
	$$
	which satisfies
	$$
	\dot{e}(t)=A(t)e(t),
	$$
	hence 
	$$
	e(t)=\Phi_A(t)\theta_e,
	$$
	with $\theta_e :=e(0)$. Consequently,
	\begin{equation*}
		\lab{x}
		x(t)=\zeta(t) - \Phi_A(t) \theta_e + G(t)\theta_\Gamma.
	\end{equation*}
	The output of the system \eqref{sys} then satisfies
	$$
	y(t)=C(\varphi(t))\left[\zeta(\varphi(t))-\Phi_A(\varphi(t))\theta_e+G(\varphi(t))\theta_\Gamma\right].
	$$
	From which we get a LRE, that allows us to identify $\theta$, as
	\begequ
	\lab{lre}
	z(t)=\Psi(t)\theta,
	\endequ
	where $\theta=\col\{\theta_e;\theta_\Gamma\}$ and the definitions from \eqref{psi} and \eqref{z} are used. It is noticeable that these signals are well defined since the assumption \eqref{boudt} ensures that $\varphi(t)$ is lower bounded. 
	
	Thirdly, we apply the DREM procedure to \eqref{lre} to generate  $n+k$ scalar LREs. Towards this end, we consider the following chain of transformations
	\begalis{
		\eqref{lre} & \Rightarrow \; \Psi^\top(t) z(t) = \Psi^\top(t)  \Psi(t) \theta \quad (\Leftarrow  \Psi^\top(t) \times) \\
		& \Rightarrow \; Y(t) = \Omega(t) \theta  \; \Big(\Leftarrow {\lambda  \over \bfp+\lambda}[\cdot]\; \mbox{and} \; \eqref{doty}, \eqref{dotome}\Big)\\
		& \Rightarrow \; \caly(t) = \Delta(t) \theta , \; (\Leftarrow \adj\{\Omega\} \times\; \mbox{and} \;  \eqref{caly}, \eqref{del} ),
	}
	with $\bfp={d \over dt}$, where we have used the fact that for any, {\em possibly singular}, square matrix $L$ we have $\adj\{L\}L=\det\{L\}I_{n+k}$ in the last line. 
	
	The last step is the analysis of the gradient estimator \eqref{dothatthe} with the FCT observer \eqref{gpebodyn}-\eqref{ftcgpebo_eta}. Replacing the latter identity in \eqref{dothatthe} provides the error equation 
	\begin{equation*}
		\lab{dottilthe}
		\dot {\tilde \theta}(t) =-\gamma \Delta^2(t) \tilde \theta(t),
	\end{equation*}
	where ${\tilde \theta}(t):= \hat \theta(t) -\theta$. The solution of this equation is given by
	\begin{equation*}
		\lab{tilthe}
		{\tilde \theta}(t)=e^{-\gamma \int_{0}^{t} \Delta^2(s)ds}\tilde \theta(0).
	\end{equation*}
	Notice that the solution of \eqref{dotw} is
	$$
	w(t)=e^{-\gamma \int_0^{t}\Delta^2(s)ds},
	$$
	hence we obtain that 
	$$
	\tilde \theta(t) =w(t)\tilde \theta(0).
	$$ 
	Clearly, this is equivalent to
	$$
	[1 - w(t)]\theta = \hat \theta(t)  - w(t) \hat \theta(0).
	$$
	On the other hand,  if $\Delta(t)$ is IE with $\rho$ satisfying \eqref{rho}, we have that there exists a $t_c \geq t_0$ such that
	$$
	w(t)=w_c(t) < 1,\; \forall t \geq t_c.
	$$
	Consequently, we conclude that 
	$$
	{1 \over 1 - w_c(t)}[\hat \theta(t) - w_c(t) \hat \theta(0)]=\theta,\; \forall t \geq t_c.
	$$
	Replacing this identity in \eqref{ftcgpebo} we complete the proof.
\end{pf}
%
\section{Discussion}
\lab{sec3}
%
Note that the upper bound on the delay implicitly impacts the value of $t_c$ given in \eqref{intexc}, thus, it predefines the fixed time of convergence.

To the best of the authors' knowledge this is the first solution to the observer problem with {\em unknown time-varying parameter} $\eta$. Observers with known  $\eta$ have been recently reported in the literature. For which the following remarks are in order.

\begite
\item	In \citep{SANGARKRS}, a more complex observer that requires strict positivity of the delay $d(t)$ and the {\em knowledge} of its derivative is reported under the classical assumption of existence of an exponentially stable Luenberger observer for the LTV system \eqref{ltvsys} with $\varphi(t)=t$, {\em i.e.} \citep[Assumption 2]{SANGARKRS}. The estimator is based on---the now classical---PDE representation of the delay, with an observer designed for the coupled LTV-PDE system. 

\item	In \citep{RUEetalijc19} a Kalman-Bucy-like observer with fractional powers is proposed, which ensures that the state estimation enters a {\em residual set} in a fixed time provided the pair $(C(t),A(t))$ is uniformly completely observable. As it is well known \citep{RUGbook}, this assumption is a {\em sufficient} condition for the verification of \citep[Assumption 2]{SANGARKRS}. It is obvious that both assumptions are {\em strictly stronger} than our IE condition.
\endite
%
\section{Simulation Results}
\lab{sec4}
%
Consider the LTV system \eqref{sys} with  $n=2$, $m=p=1$, $q=2$, $k=3$ and
\begalis{
	&\cala(y,t) =
	\begmat{
		-y^2& 1\\
		-\sin^2(t) & 0
	},\; 
\calb(u,y)=
	\begmat{
		0\\y^3u
	},\; \\
&C=\begmat{
		1\\0
	}, \;
\cald(y,t)=\begmat{
	-2\sin(t) &0\\0 &-y^3
},
\eta(t)=\begmat{
	1+0.1 \sin(3t)\\1+0.3 \cos(3t)
}.
}
It is obviously vector $\eta$ can be generated in the form \eqref{eta}, \eqref{xi} with 
	$H(t)=
	\begmat{
		&1 \; &0 \; &1\\&0 &1 &1
	}$, 
$\Gamma(t)=
\begmat{
&0 \; &1 \; &0\\&-9 &0 &0\\&0 &0 &0
},$ and vector of initial condition $\xi(0)=\col(0 \; 0.3 \; 1)$

We consider the input $u=-1$, the state initial conditions $x(0)=\col(1,2)$ and zero for all the other variables. For the simulations we used $\lambda=1$, $\mu=0.01$ and different values of the adaptation gain $\gamma$. Selection of parameters $\lambda$, $\gamma$ and $\mu$  is a complicated problem, we need to find such values of parameters  which ensure boundedness of all signals  and convergence of estimation errors to zero. Convergence time $t_c$ depends on the values of $\mu$ and $\gamma$.

\begin{figure}[ht]
	\centering
	\includegraphics[width=1.1\linewidth]{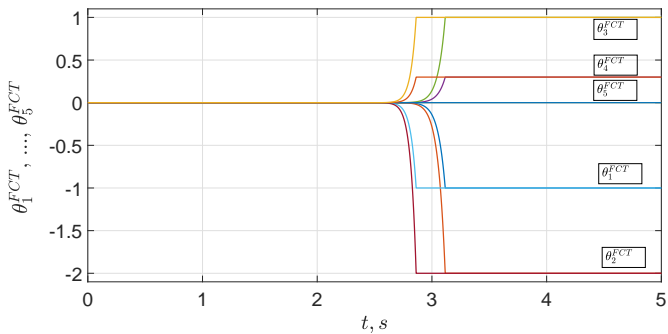}
	\caption{Transients of unknown parameter estimations ${\theta}^{FCT}$ for  $\gamma=10^{10}$ and $\gamma=10^{12}$ and case {\bf C1}}
	\label{fig_0}
\end{figure}

\begin{figure}[ht]
	\centering
	\includegraphics[width=1.1\linewidth]{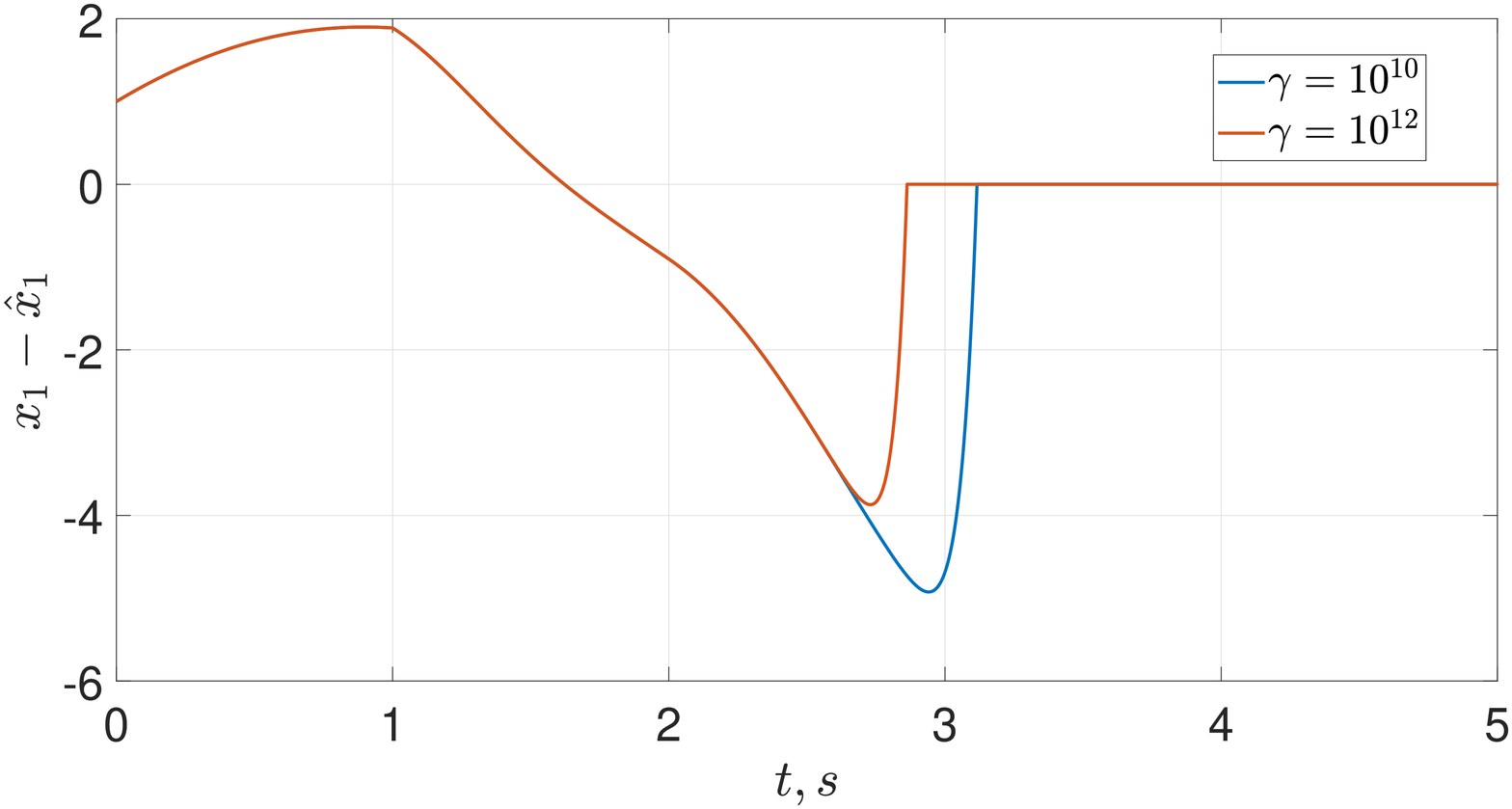}
	\caption{Error transients $x_1(t)-\hat{x}_1(t)$ for different $\gamma$ and case {\bf C1}}
	\label{fig_1}
\end{figure}

\begin{figure}[ht]
	\centering
	\includegraphics[width=1.1\linewidth]{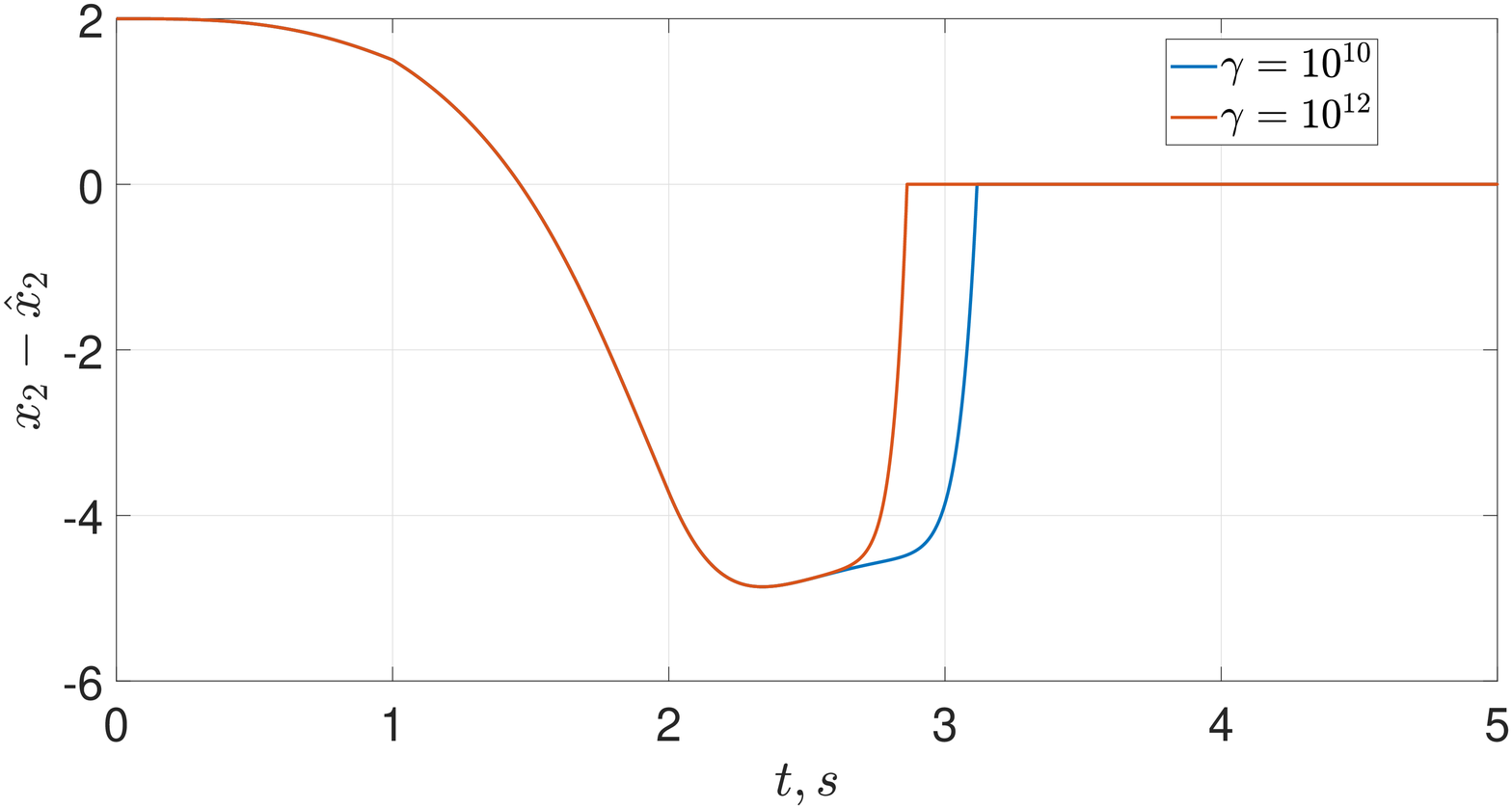}
	\caption{Error transients $x_2(t)-\hat{x}_2(t)$ for different $\gamma$ and case {\bf C1}}
	\label{fig_2}
\end{figure}

\begin{figure}[ht]
	\centering
	\includegraphics[width=1.1\linewidth]{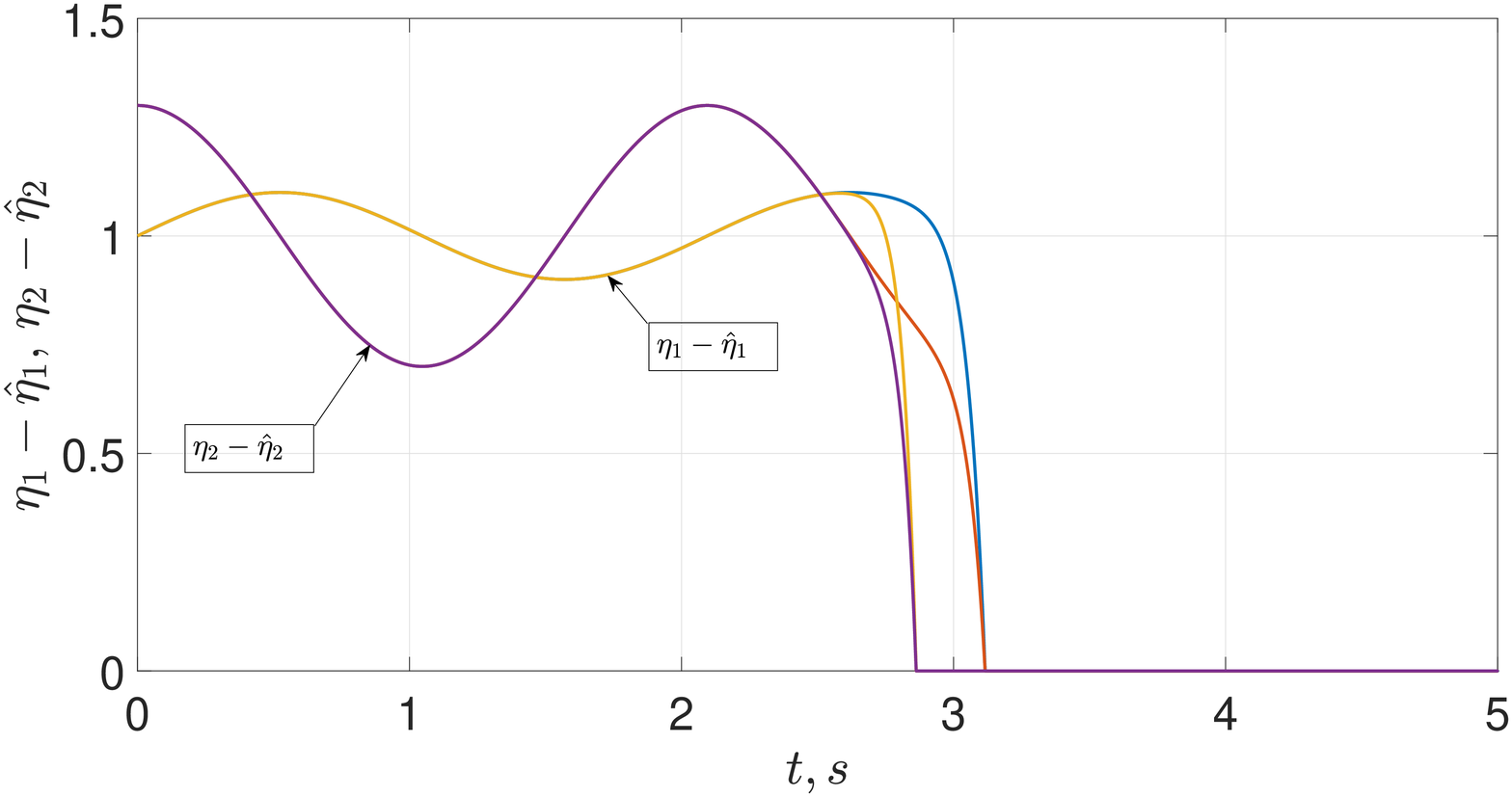}
	\caption{Errors transients $\eta_1(t)-\hat{\eta}_1(t)$ and $\eta_2(t)-\hat{\eta}_2(t)$ for $\gamma=10^{10}$ and $\gamma=10^{12}$ and case {\bf C1}}
	\label{err_eta_C1}
\end{figure}

\begin{figure}[ht]
	\centering
	\includegraphics[width=1.1\linewidth]{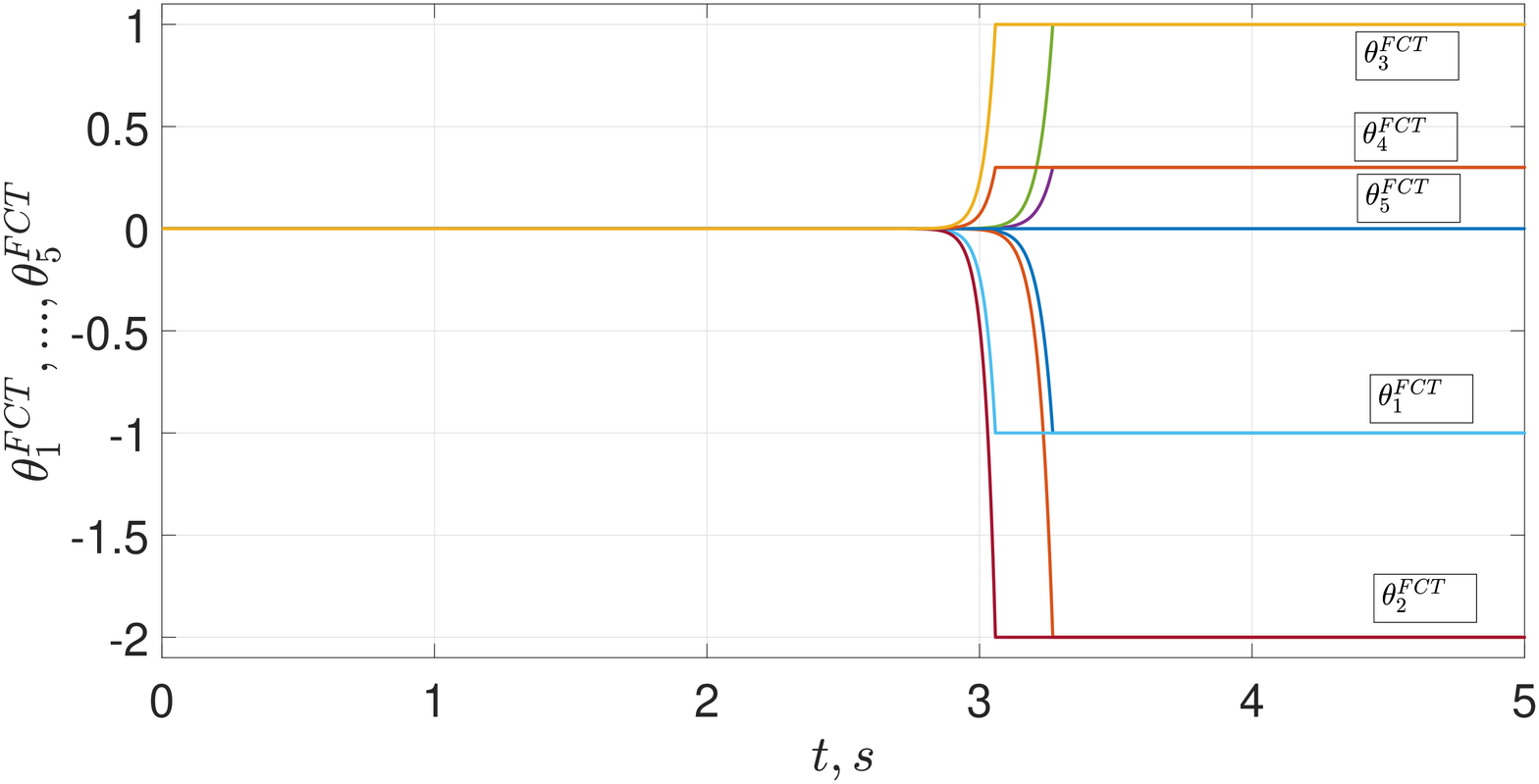}
	\caption{Transients of unknown parameter estimations ${\theta}^{FCT}$ for  $\gamma=10^{10}$ and $\gamma=10^{12}$ and case {\bf C2}}
	\label{fig_02}
\end{figure}

\begin{figure}[ht]
	\centering
	\includegraphics[width=1.1\linewidth]{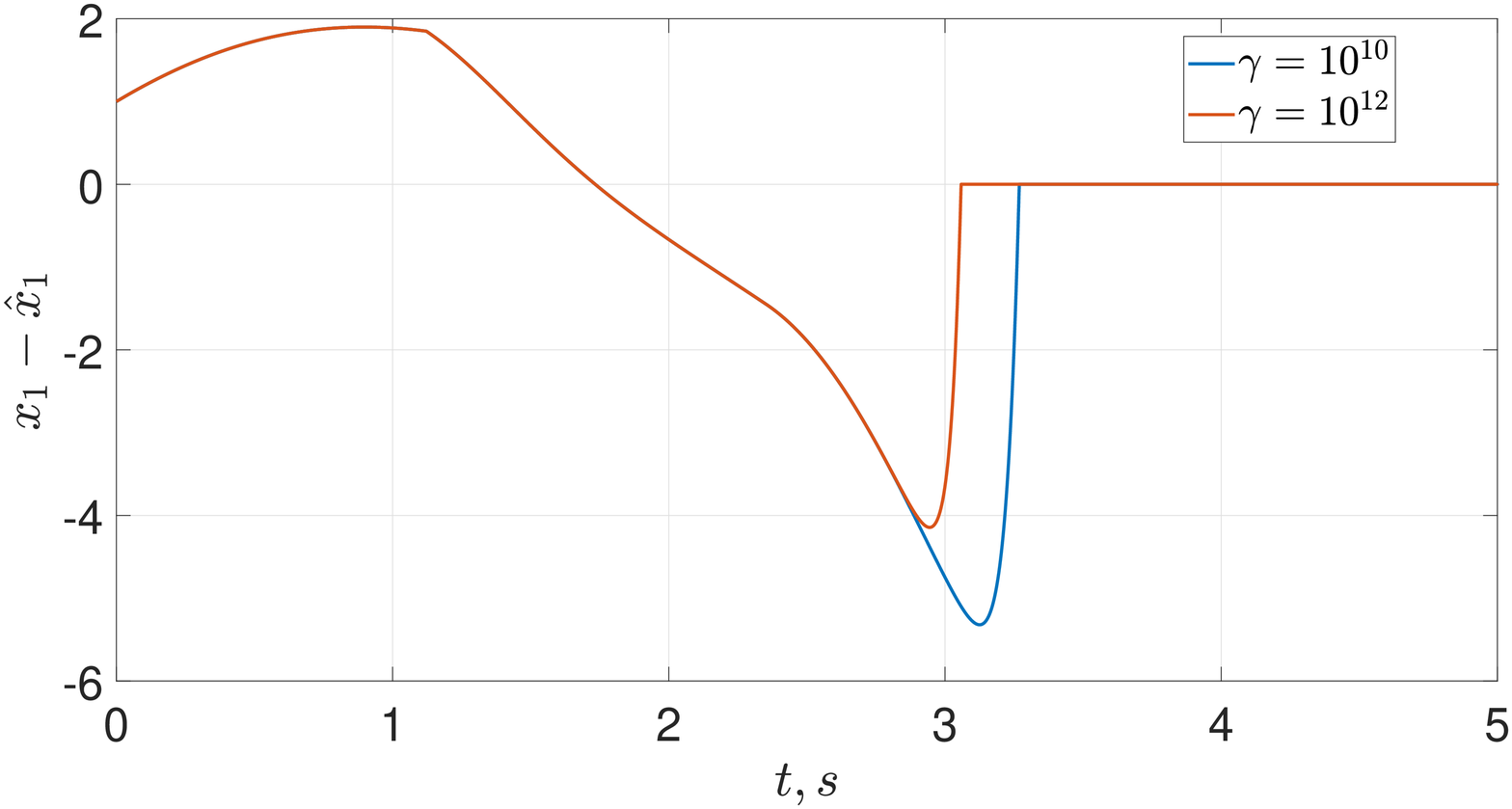}
	\caption{Error transients $x_1(t)-\hat{x}_1(t)$ for different $\gamma$ and case {\bf C2}}
	\label{fig_3}
\end{figure}

\begin{figure}[ht]
	\centering
	\includegraphics[width=1.1\linewidth]{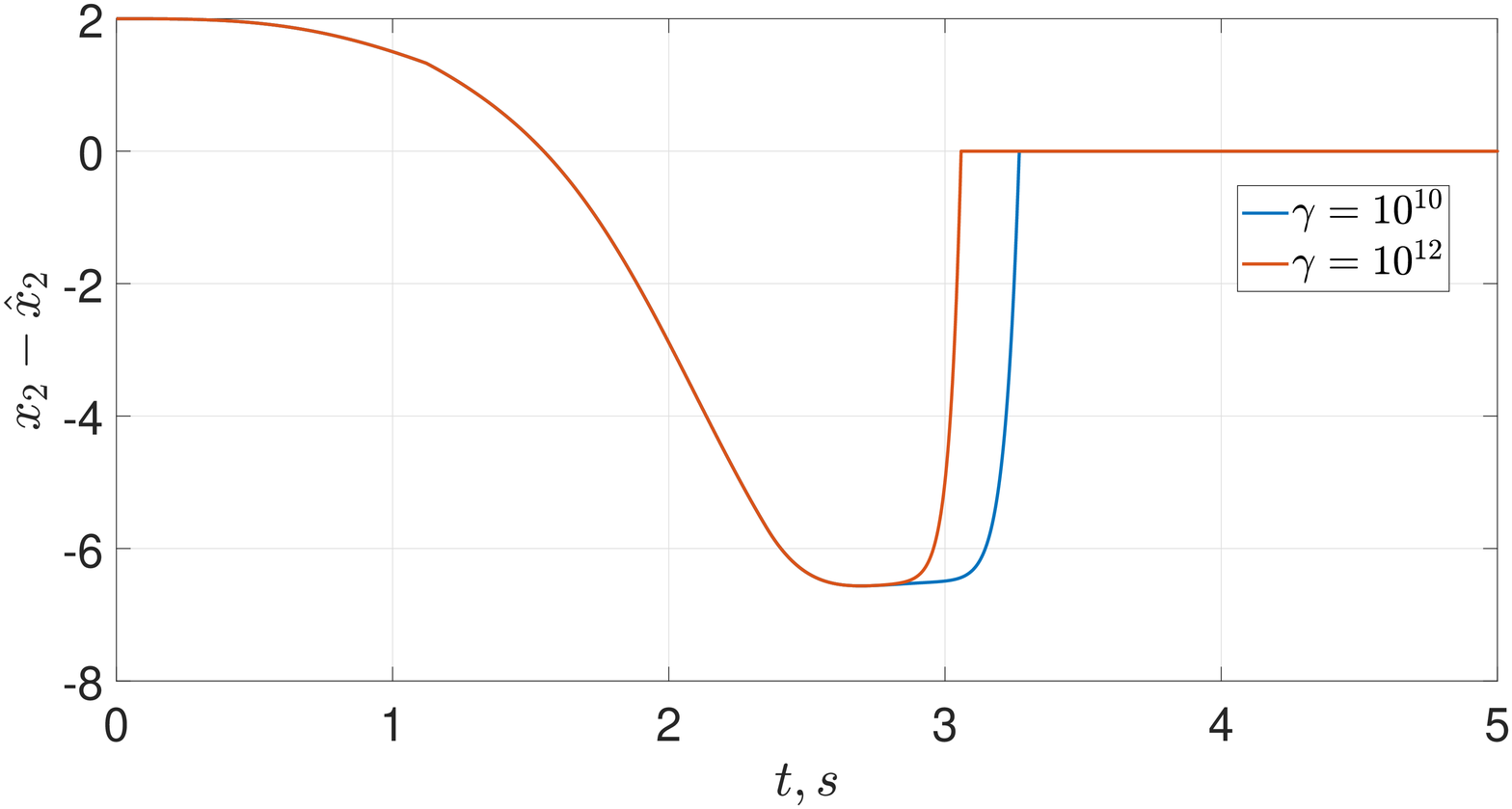}
	\caption{Error transients $x_2(t)-\hat{x}_2(t)$ for diffrerent $\gamma$ and case {\bf C2}}
	\label{fig_4}
\end{figure}

\begin{figure}[ht]
	\centering
	\includegraphics[width=1.1\linewidth]{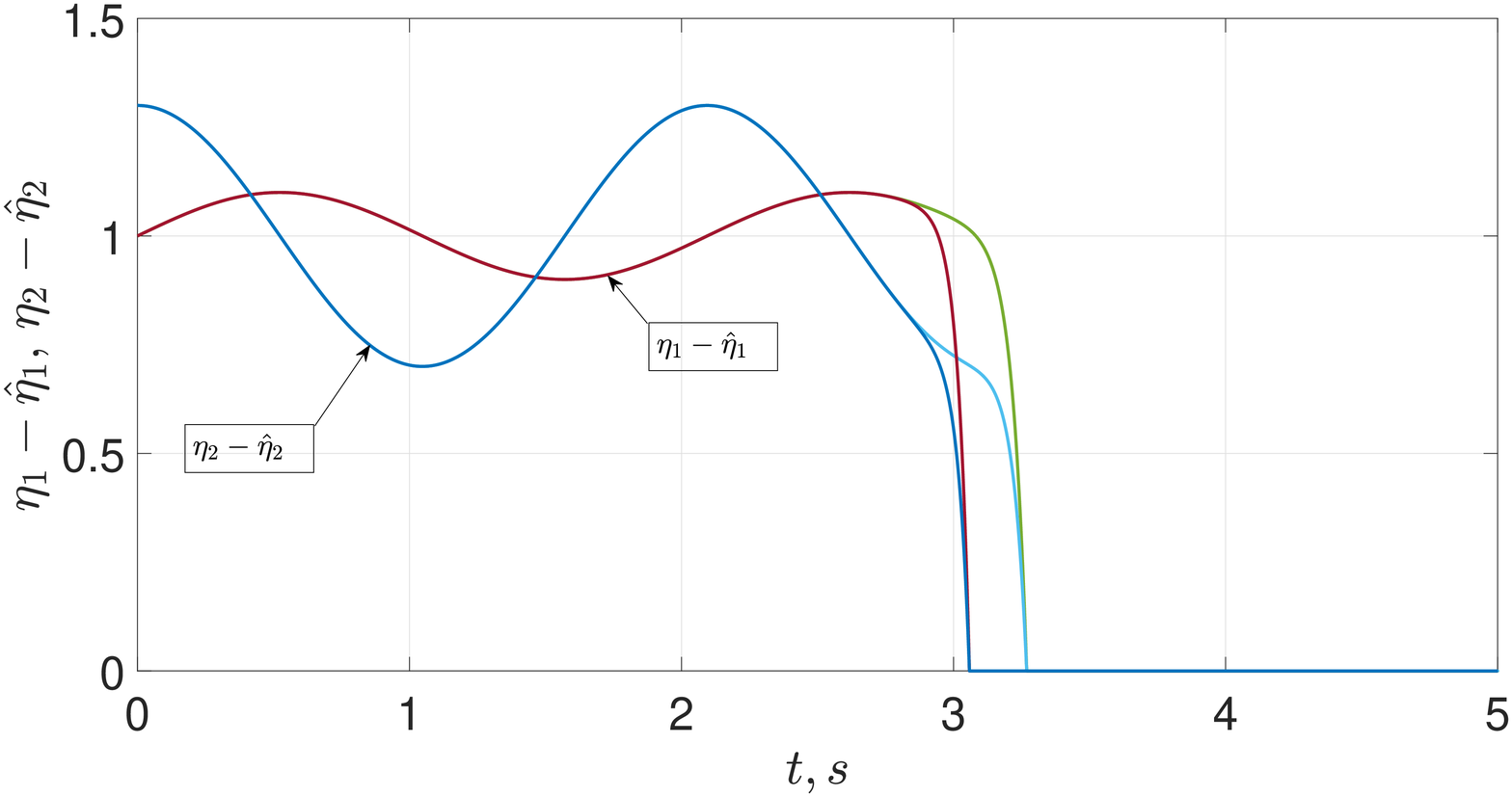}
	\caption{Errors transients $\eta_1(t)-\hat{\eta}_1(t)$ and $\eta_2(t)-\hat{\eta}_2(t)$ for $\gamma=10^{10}$ and $\gamma=10^{12}$ and case {\bf C2}}
	\label{err_eta_C2}
\end{figure}

\begin{figure}[ht]
	\centering
	\includegraphics[width=1.1\linewidth]{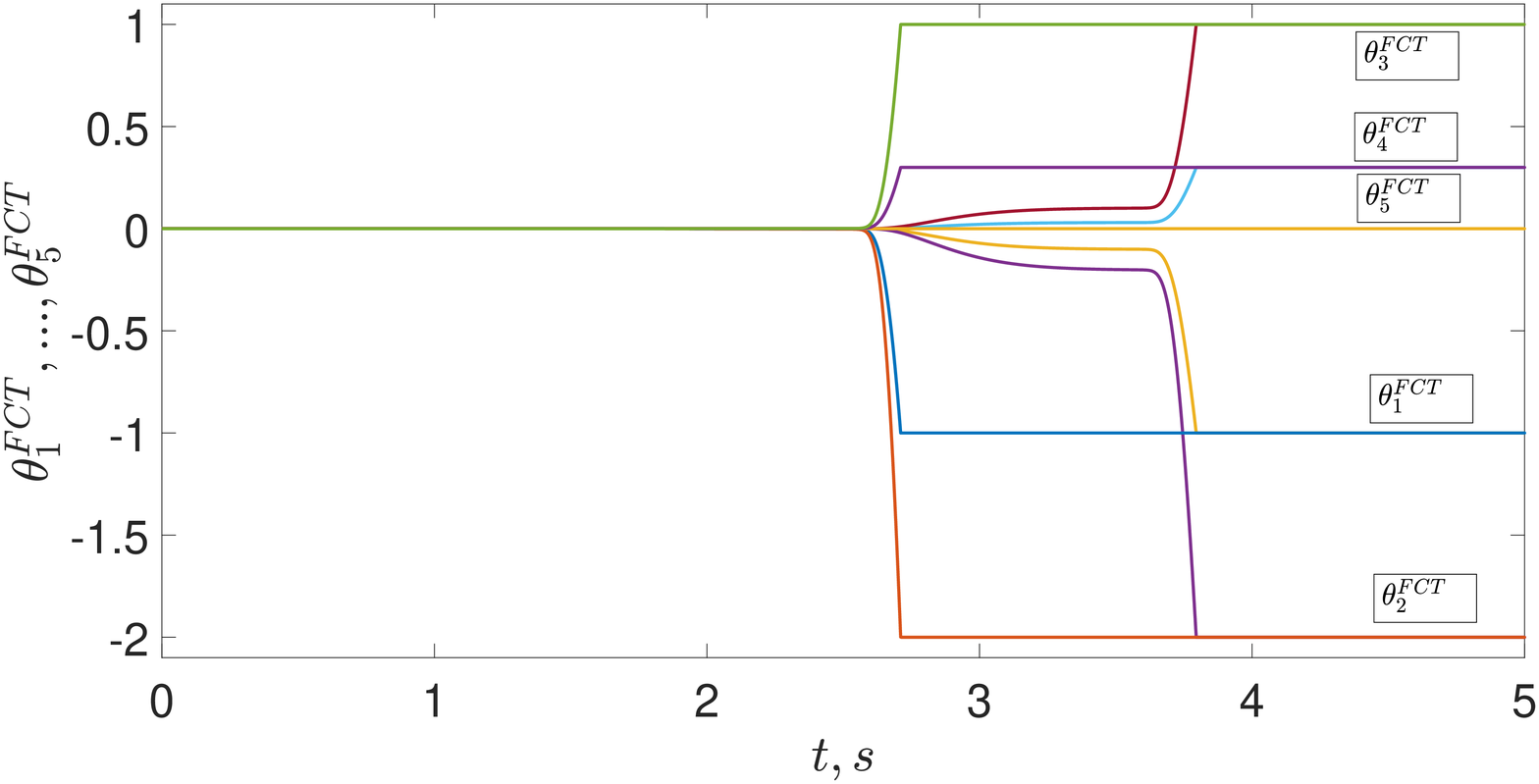}
	\caption{Transients of unknown parameter estimations ${\theta}^{FCT}$ for  $\gamma=10^{10}$ and $\gamma=10^{12}$ and case {\bf C3}}
	\label{fig_03}
\end{figure}

\begin{figure}[ht]
	\centering
	\includegraphics[width=1.1\linewidth]{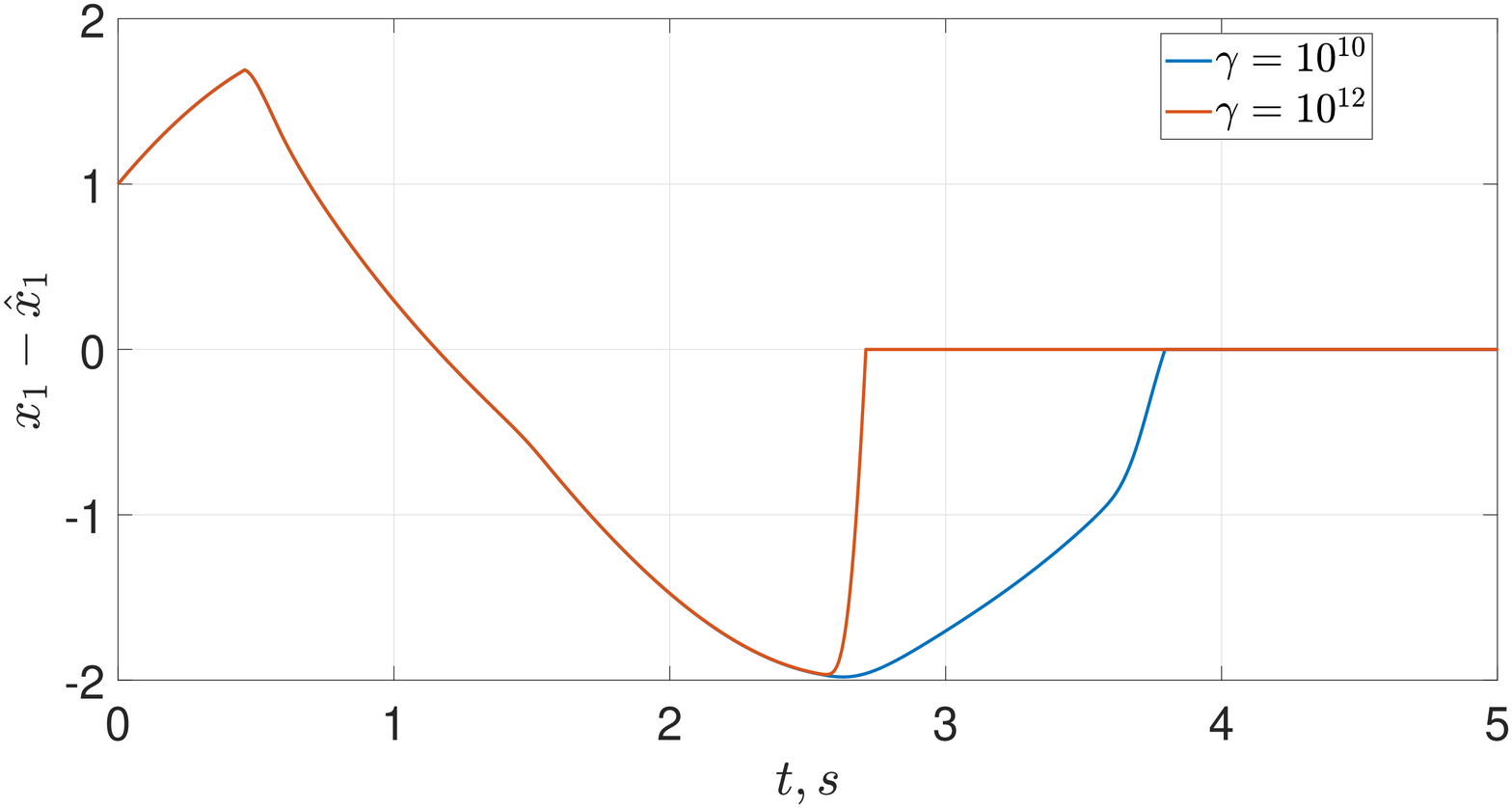}
	\caption{Error transients $x_1(t)-\hat{x}_1(t)$ for diffrerent $\gamma$ and case {\bf C3}}
	\label{fig_5}
\end{figure}

\begin{figure}[ht]
	\centering
	\includegraphics[width=1.1\linewidth]{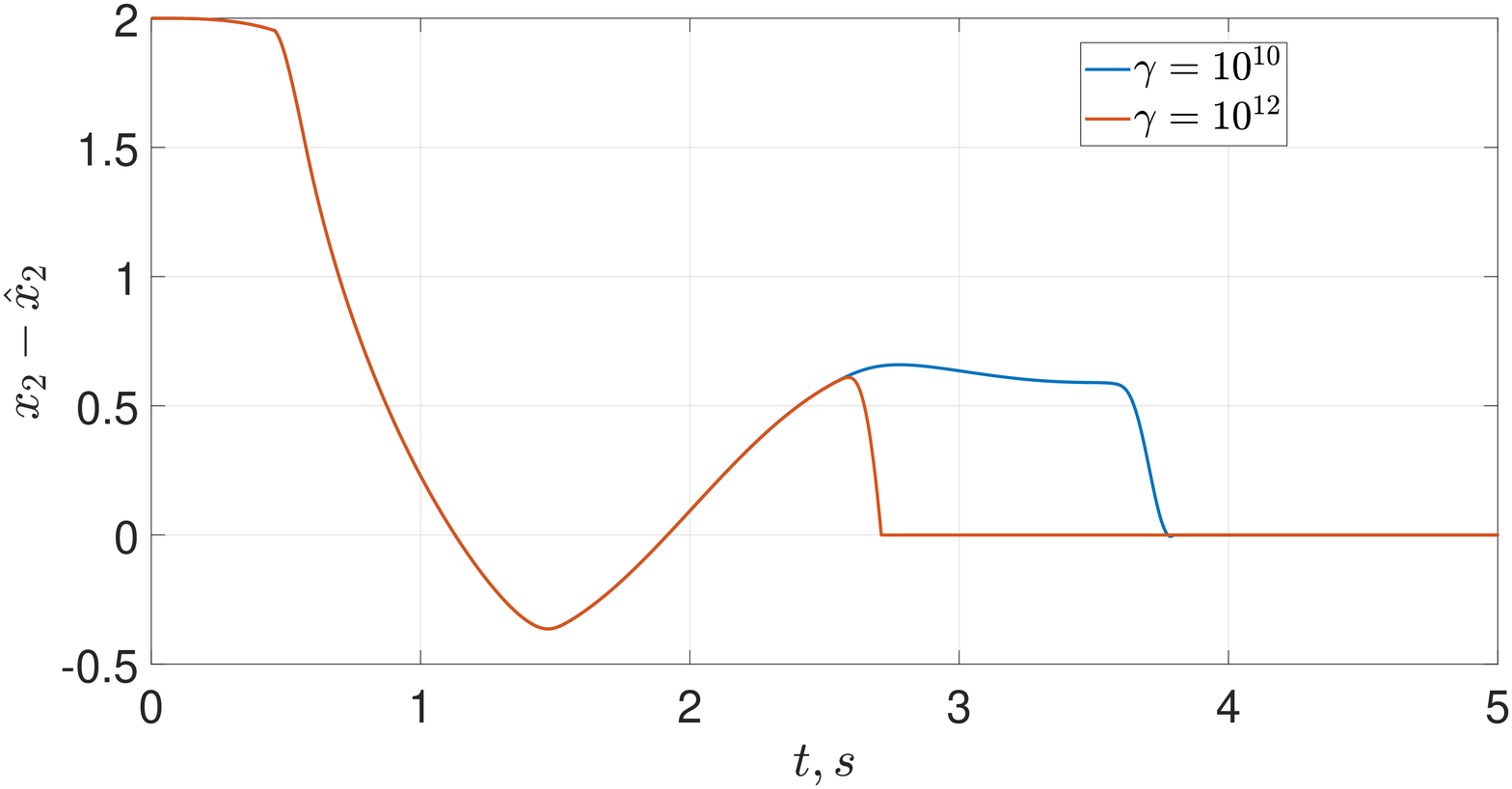}
	\caption{Error transients $x_2(t)-\hat{x}_2(t)$ for diffrerent $\gamma$ and case {\bf C3}}
	\label{fig_6}
\end{figure}

\begin{figure}[ht]
	\centering
	\includegraphics[width=1.1\linewidth]{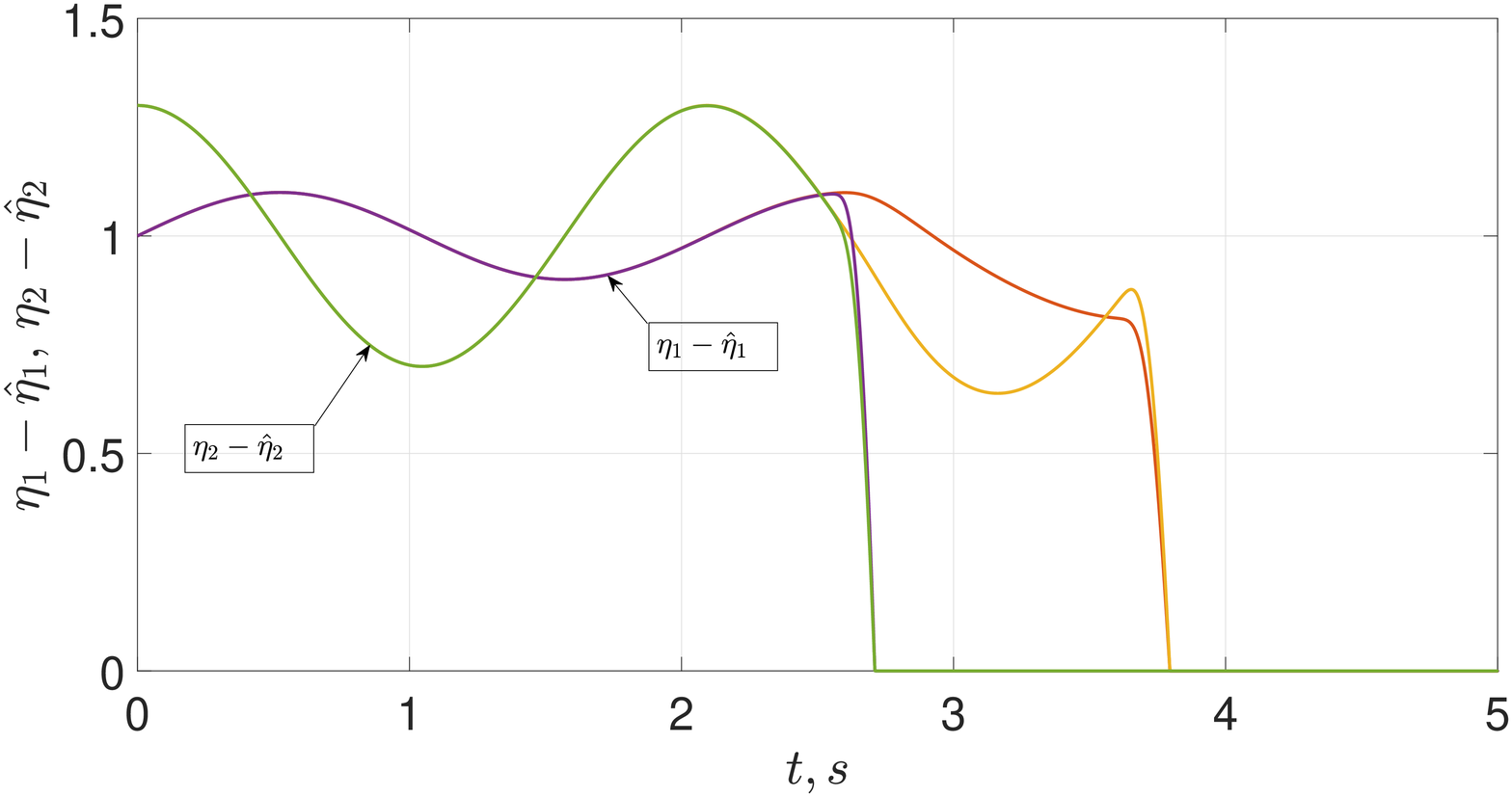}
	\caption{Errors transients $\eta_1(t)-\hat{\eta}_1(t)$ and $\eta_2(t)-\hat{\eta}_2(t)$ for $\gamma=10^{10}$ and $\gamma=10^{12}$ and case {\bf C3}}
	\label{err_eta_C3}
\end{figure}

We consider three cases for the delay function \eqref{boudt}:
\begenu[{\bf C1}]
\item $d(t)=1$ (Fig. \ref{fig_0}, Fig. \ref{fig_1}, Fig. \ref{fig_2} and Fig. \ref{err_eta_C1});
\item $d(t)=1+0.25\sin(t)$ (Fig. \ref{fig_02}, Fig. \ref{fig_3}, Fig. \ref{fig_4} and Fig. \ref{err_eta_C2});
\item $d(t)=0.1+\cos^2(3t)$ (Fig. \ref{fig_03}, Fig. \ref{fig_5}, Fig. \ref{fig_6} and Fig. \ref{err_eta_C3}).
\endenu

From simulation results shown in Fig. \ref{fig_0}, Fig. \ref{fig_02}, Fig. \ref{fig_03} we can see that the estimations of parameters converge to the values of the vector of initial conditions of unknown time-varying parameter $\eta=\col(0, 0.3, 1)$ and the state initial conditions $x(0)=\col(1,2)$ in fixed-time, which can be reduced by increasing the value of adaptation gain $\gamma$. From simulation results shown in Fig. \ref{fig_1}, Fig. \ref{fig_2},Fig. \ref{err_eta_C1}, Fig. \ref{fig_3}, Fig. \ref{fig_4}, Fig. \ref{err_eta_C2}, Fig. \ref{fig_5}, Fig. \ref{fig_6}, Fig. \ref{err_eta_C3} we can see that the estimation errors of the state variables $x$ and $\eta$ also converge to zero in fixed time.


%
\section{Conclusions and Future Research}
\lab{sec5}
%
An adaptive state observer for time-varying affine-in-the-state systems with unknown time-varying paramters and  delayed measurements of the form \eqref{sys} has been proposed. Due to the use of DREM, this observer ensures fixed-time convergence under an extremely weak assumption of IE. Following the PEBO approach, the observation of the state is carried out via the estimation of some suitable initial conditions. To the best of our knowledge this is the only result available for this problem.

Current research is under way to relax the strict assumption of exact knowledge of the time delay function. Robustness properties of the proposed estimator with respect to measurement noise and plant disturbances can be evaluated.

\bibliographystyle{plain}        
\bibliography{alcos}             

\end{document}